\documentstyle[12pt,amssymb,amsbsy,righttag]{amsart}
\newcommand{\g}{\gamma}

\newcommand{\z}{\zeta}

\newcommand{\r}{\rho}

\newcommand{\f}{\varphi}
\renewcommand{\o}{\omega}

\renewcommand{\O}{\Omega}
\newcommand{\U}{\Upsilon}

\newcommand{\A}{{\cal A}}

\newcommand{\I}{{\cal I}}

\newcommand{\cL}{{\cal L}}

\newcommand{\1}{\boldsymbol{1}}

\newcommand{\C}{{\Bbb C}}
\newcommand{\T}{{\Bbb T}}
\newcommand{\pp}{{\Bbb P}}
\newcommand{\dd}{{\Bbb D}}
\newcommand{\R}{{\Bbb R}}
\newcommand{\Z}{{\Bbb Z}}
\newcommand{\mm}{{\Bbb M}}
\newcommand{\0}{{\Bbb O}}

\newcommand{\rf}[1]{(\ref{#1})}

\newcommand{\df}{\stackrel{\mathrm{def}}{=}}
\newcommand{\dist}{\operatorname{dist}}
\newcommand{\Ker}{\operatorname{Ker}}

\newcommand{\rank}{\operatorname{rank}}

\newcommand{\eeq}{\end{equation}}
\newcommand{\beq}{\begin{equation}}
\newcommand{\bay}{\begin{eqnarray}}
\newcommand{\ey}{\end{eqnarray}}

\newcommand{\be}{\infty}

\newcommand{\bl}{\blacksquare}
\newcommand{\ess}{\operatorname{ess}}
\newcommand{\ind}{\operatorname{ind}}
\newcommand{\wind}{\operatorname{wind}}

\newcommand{\Pf}{{\bf Proof. }}

\newcommand{\ov}{\overline}

\newtheorem{thm}{\hspace{\parindent}Theorem}[section]

\newtheorem{cor}[thm]{\hspace{\parindent}Corollary}
\newtheorem{lem}[thm]{\hspace{\parindent}Lemma}

\begin{document}
\newcommand{\bs}{\boldsymbol}
\newcommand{\vse}{\vspace{.2in}}
\renewcommand{\theequation}{\thesection.\arabic{equation}}

\author{R.B. Alexeev and V.V. Peller}
\thanks{The second author is partially supported by NSF grant DMS 9970561}

\title{Invariance properties of thematic\\ factorizations of matrix functions}

\maketitle

%



\begin{abstract}
We study the problem of invariance of indices of thematic factorizations.
Such factorizations were introduced in [PY1] for studying superoptimal
approximation by bounded analytic matrix functions. As shown in [PY1], the indices
may depend on the choice of a thematic factorization. We introduce the notion
of a monotone thematic factorization. The main result shows that under natural
assumptions a matrix function that admits a thematic factorization also admits 
a monotone thematic factorization and the indices of a monotone thematic 
factorization are uniquely determined by the matrix function itself. We obtain 
similar results for so-called partial thematic factorizations.
\end{abstract}

\setcounter{equation}{0}
\section{\bf Introduction}

\

It is well known [Kh] that for a continuous scalar function $\f$ on the
unit circle $\T$ there exists a unique function $f\in H^\be$ such that
$$
\dist_{L^\be}(\f,H^\be)=\|\f-f\|_{L^\be}.
$$
However, the situation in the case of matrix-valued function is considerably
more complicated.

Suppose that $\Phi$ is a matrix function in $L^\be(\mm_{m,n})$, i.e., $\Phi$ is
an essentially bounded function on the unit circle $\T$ 
that takes values in the space $\mm_{m,n}$
of $m\times n$ matrices. We say that a function $F\in H^\be(\mm_{m,n})$ (by this we 
mean that all entries of $F$ belong to $H^\be$) is a {\it best approximation of}
$\Phi$ by bounded analytic matrix functions if
$$
\|\Phi-F\|_{L^\be}=\dist_{L^\be}(\Phi,H^\be(\mm_{m,n})).
$$
Here for a function $\Psi$ in $L^\be(\mm_{m,n})$ we use the notation
$$
\|\Psi\|_{L^\be}\df\ess\sup_{\z\in\T}\|\Psi(\z)\|_{\mm_{m,n}},
$$
where $\mm_{m,n}$ is equipped with the operator norm from $\C^n$ to $\C^m$.

It is easy to see that unlike the scalar case we can have uniqueness only
in exceptional cases. Indeed, if
$\Phi=\left(\begin{array}{cc}\bar z&0\\0&0\end{array}\right)$, then
$\dist_{L^\be}(\Phi,H^\be(\mm_{2,2}))=1$ since $\dist_{L^\be}(\bar z,H^\be)=1$. 
Clearly, for any scalar function $f\in H^\be$ such that $\|f\|_\be\le1$ we have
$$
\left\|\left(\begin{array}{cc}\bar z&0\\0&-f\end{array}\right)\right\|_{L^\be}=1,
$$
and so $\left(\begin{array}{cc}0&0\\0&f\end{array}\right)$ is a best approximation
of $\Phi$.

Recall that by a matrix analog of Nehari's theorem (see [Pa]),
$$
\dist_{L^\be}(\Phi,H^\be(\mm_{m,n}))=\|H_\Phi\|,
$$
where the {\it Hankel operator} $H_\Phi:H^2(\C^n)\to H^2_-(\C^m)$
is defined by
$$
H_\Phi f\df\pp_-\Phi f,\quad f\in H^2(\C^n).
$$
Here $\pp_-$ is the orthogonal projection onto 
$H^2_-(\C^m)\df L^2(\C^m)\ominus H^2(\C^m)$.

Recall also that by Hartman's theorem (see e.g., [N]), 
$H_\Phi$ is compact if and only if
$\Phi\in(H^\be+C)(\mm_{m,n})$, where
$$
H^\be+C\df\{f+g:~f\in H^\be,~g\in C(\T)\}.
$$
(Throughout the paper we write $\Phi\in X(\mm_{m,n})$ if all entries of an 
$m\times n$ matrix function $\Phi$ belong to a function space $X$; sometimes
to simplify the notation we will write simply $\Phi\in X$ if this does not
lead to a confusion.)

In [PY1] it was shown that if $\Phi\in(H^\be+C)(\mm_{m,n})$, then there exists
a unique function $F\in H^\be(\mm_{m,n})$ that minimizes (lexicographically)
not only $\|\Phi-F\|_{L^\be}$ but also the essential {\it suprema}
$$
t_j\df\ess\sup_{\z\in\T}s_j(\Phi(\z)-F(\z)),\quad j\le\min\{m,n\}-1
$$
of all subsequent singular values of $\Phi(\z)-F(\z)$, $\z\in\T$. Such functions
$F$ are called {\it superoptimal approximations of} $\Phi$ by bounded analytic
matrix functions. The numbers $t_j$ are called {\it superoptimal singular
values of} $\Phi$. It was also shown in [PY1] that the error function
$\Phi-F$ admits certain special factorizations ({\it thematic factorizations}).
For each such factorization the sequence of positive indices $k_j$, \linebreak
$j\ge0,\,t_j>0$, ({\it thematic indices}) was defined. We refer the reader to \S 2 
where formal definitions are given. Note that another approach to superoptimal 
approximation was found later in [T].

In [PT2] the same results were proved for functions $\Phi\in L^\be(\mm_{m,n})$
such that the essential norm $\|H_\Phi\|_{\text e}$ of $H_\Phi$ (i.e., the
distance from $H_\Phi$ to the set of compact operators) is less than the smallest
nonzero superoptimal singular value of $\Phi$. Recall that
$$
\|H_\Phi\|_{\text e}=\dist_{L^\be}\big(\Phi,(H^\be+C)(\mm_{m,n})\big)
$$
(see e.g., [S] for the proof of this formula for scalar functions, in the
matrix-valued case the proof is the same).

It turned out, however, that the thematic indices are not uniquely determined
by the function $\Phi$ itself but may depend on the choice of a thematic 
factorization (see [PY1]). On the other hand it was shown in [PY2] (see also [PT2]) 
that the sum of the thematic indices that correspond to the superoptimal singular
values equal to a specific number is uniquely determined by $\Phi$.

In this paper we show that one can always choose a so-called monotone thematic 
factorization, i.e., a thematic factorization such that the indices that correspond
to equal superoptimal nonzero singular values are arranged in the nonincreasing 
order. We refer the reader to \S 4 for a formal definition. We prove in \S 3 and
\S 4 that the indices of a monotone thematic factorization are uniquely determined 
by the function $\Phi$ itself. Section 2 contains
definitions and statements of basic results on superoptimal approximation and 
thematic factorizations.

Note that using the same methods we can obtain similar results in the case
of the four block problem (which is an important generalization of the 
problem of best approximation by bounded analytic matrix functions). We refer
the reader to [PT2] which contains results on superoptimal approximation
and thematic factorizations related to the four block problem.

We can also obtain similar results in the case of infinite matrix functions.
We refer the reader to [T], [Pe], and [PT1] for results on superoptimal approximation
and thematic factorizations for infinite matrix functions.

\

\setcounter{equation}{0}
\section{\bf Superoptimal approximation and thematic factorizations}

\

In this section we collect necessary information on superoptimal approximation
and thematic factorizations.

Let $\Phi\in L^\be(\mm_{m,n})$. We put
$$
\O_0=\{F\in H^\be(\mm_{m,n})
:~F~\mbox{minimizes}~ t_0\df\ess\sup_{\z\in\T}\|\Phi(\z)-F(\z)\|\};
$$
$$
\O_j=\{F\in \O_{j-1}:~F~\mbox{minimizes}~ 
t_j\df\ess\sup_{\z\in\T}s_j(\Phi(\z)-F(\z))\}.
$$
Recall that for a matrix $A\in\mm_{m,n}$ the {\it$j$th singular value} $s_j(A)$
is defined by
$$
s_j(A)=\inf\{\|A-R\|:~\rank R\le j\},\quad j\ge0.
$$ 

Functions $F$ in $\O_{\min\{m,n\}-1}$ are called {\it superoptimal approximations
of $\Phi$ by analytic functions}, or superoptimal solutions of the Nehari problem.
The numbers $t_j$ are called {\it superoptimal singular values of} $\Phi$.
The notion of superoptimal approximation plays an important role
in $H^\be$ control theory. 

It can be shown easily with the help of a compactness argument that the sets
$\O_j$ are nonempty. In particular, for any matrix function in $L^\be(\mm_{m,n})$
there exists a superoptimal approximation by analytic matrix functions.

It was shown in [PY1] that for any matrix function $\Phi\in(H^\be+C)(\mm_{m,n})$
there exists a unique superoptimal approximation.
We denote by $\A\Phi$ the unique superoptimal 
approximation of $\Phi$ by bounded analytic matrix functions whenever it is unique.

Later in [PT2] stronger results were obtained. It was shown there that if \linebreak
$\Phi\in L^\be(\mm_{m,n})$ and the essential norm $\|H_\Phi\|_{\text e}$
of the Hankel operator $H_\Phi$ is less than the smallest nonzero superoptimal
singular value of $\Phi$, then $\Phi$ has a unique superoptimal approximation
by bounded analytic matrix functions.

A matrix function $\Phi\in L^\be(\mm_{m,n})$ is called {\it badly approximable}
if
$$
\dist_{L^\be}(\Phi,H^\be(\mm_{m,n}))=\|\Phi\|_{L^\be}.
$$
It is called {\it very badly approximable} if the zero matrix function is a 
superoptimal approximation of $\Phi$.

Recall that a nonzero scalar function $\f\in H^\be+C$ is badly approximable if and 
only if it has constant modulus almost everywhere on $\T$, belongs to $QC$, and its
winding number $\wind\f$ is negative, where the space $QC$ of quasi-continuous
functions is defined by
$$
QC=\{f\in H^\be+C:~\bar f\in H^\be+C\}.
$$ 
For continuous $\f$ this was proved in [Po]
(see also [AAK1]). For the general case see [PK]. Recall that if $\f\in QC$ and
$\f$ has constant modulus on $\T$ almost everywhere, the harmonic extension
of $\f$ to the unit disk $\dd$ is separated away from zero near the unit circle
and $\wind\f$ is defined as the winding number of the restriction of the harmonic
extension of $\f$ to the circle of radius $\r$ for $\r$ sufficiently close to 1 
(see [D]). Note also that if $\f\in QC$ and $\f$ has constant modulus on $\T$, then
the Toeplitz operator $T_\f$ on $H^2$ is Fredholm and its index $\ind T_\f$
equals $-\wind\f$ (see [D]). Recall that for $\f\in L^\be$ the Toeplitz operator
$T_\f$ on $H^2$ is defined by
$$
T_\f f=\pp_+\f f,\quad f\in H^2,
$$
where $\pp_+$ is the orthogonal projection onto $H^2$.

A similar description holds for functions $\f\in L^\be$ such that
$\|H_\f\|_{\text e}<\|H_\f\|$. In this case $\f$ is badly approximable if
and only if $\f$ has constant modulus almost everywhere on $\T$, the Toeplitz
operator $T_\f$ is Fredholm and $\ind T_\f>0$.

To state the description of badly approximable and very badly approximable 
matrix functions obtained in [PY1] and [PT2], we need the notion of a thematic 
matrix function. Recall that a function $F\in H^\be(\mm_{m,n})$ is called {\it inner}
if \linebreak$F^*(\z)F(\z)=I_n$ almost everywhere on $\T$ ($I_n$ stands for the 
identity matrix in $\mm_{n,n}$). $F$ is called {\it outer} if $FH^2(\C^n)$ is dense 
in $H^2(\C^m)$. Finally, $F$ is called {\it co-outer} if the transposed function 
$F^{\text t}$ is outer.

An $n\times n$ matrix function $V$, $n\ge2$, is called {\it thematic} if it is 
unitary-valued and has the form
$$
V=\left(\begin{array}{cc}\bs{v}&\ov{\Theta}\end{array}\right),
$$
where the matrix functions $\bs{v}\in H^\be(\C^n)$ and $\Theta\in H^\be(\mm_{n,n-1})$
are both inner and co-outer. Note that if $V$ is a thematic function, then
all minors of $V$ on the first column (i.e., minors of an arbitrary size
that involve the first column) belong to $H^\be$ ([PY1]). If $n=1$, a 
thematic function is a constant function whose modulus is equal to 1.

It was shown in [PY1] that a function 
$\Phi\in(H^\be+C)(\mm_{m,n})\setminus H^\be(\mm_{m,n})$ 
is badly approximable if and only if it admits a representation
\beq
\label{2.1}
\Phi=W^*\left(\begin{array}{cc}su&0\\0&\Psi\end{array}\right)V^*,
\end{equation}
where $s>0$, $V$ and $W^{\text t}$ are thematic functions, $u$ is a scalar 
{\it unimodular} function (i.e., $|u(\z)|=1$ for almost all $\z\in\T$) in $QC$ with 
negative winding number, and \linebreak$\|\Psi\|_{L^\be}\le s$. Note that in this 
case $V$ and $W$ must belong to $QC$, $\Psi$ must belong to $H^\be+C$, and 
$s=\|H_\Phi\|$ (see [PY1]).

A similar result was obtained in [PT2] in the more general case when \linebreak
$\|H_\Phi\|_{\text e}<\|H_\Phi\|$. Such a matrix function $\Phi$ is badly 
approximable if and only if it admits a representation of the form (\ref{2.1})
in which $s>0$, $\|\Psi\|_{L^\be}\le s$, $V$ and $W^{\text t}$ are
thematic matrix functions , and $u$ is a unimodular function such that $T_u$ is 
Fredholm and $\ind T_u>0$.

Suppose now that $m\le n$. It was proved in [PY1] that a matrix function \linebreak
$\Phi\in(H^\be+C)(\mm_{m,n})$ is very badly approximable if and only if
$\Phi$ admits a representation
\beq
\label{2.2}
\!\!\!\!\!\!
\Phi = W_0^*\cdots W^*_{m-1} 
\left( \begin{array}{ccccccc}
s_0 u_0 & 0 & \cdots & 0 & 0 & \cdots & 0 \\
0 & s_1 u_1 & \cdots & 0 & 0 & \cdots & 0\\
\vdots & \vdots & \ddots & \vdots & \vdots & \ddots & \vdots\\
0 & 0 & \cdots & s_{m-1} u_{m-1} & 0 & \cdots & 0 \end{array}
\right)
V^*_{m-1}\cdots V^*_0
\end{equation}
for some badly approximable unimodular functions $u_0, \cdots, u_{m-1} \in QC$
and some nonincreasing sequence $\{s_j\}_{0\le j\le m-1}$ of nonnegative numbers;
\beq
\label{2.3}
W_j=\left(\begin{array}{cc}I_j & 0 \\ 0 & \breve{W}_j \end{array} \right),\quad
V_j = \left( \begin{array}{cc} I_j & 0 \\ 0 & \breve{V}_j \end{array} \right),
\quad 1 \leq j \leq m-1,
\end{equation}
and $W^{\text t}_0, \breve{W}^{\text t}_j, V_0, \breve{V}_j$ are thematic 
matrix functions, 
$1 \le j\le m-1$. Moreover, in this case the $s_j$ are the superoptimal
singular values of $\Phi$: $s_j=t_j$, $0\le j\le m-1$, and the 
matrix functions $V_j,\,W_j$, $0\le j\le m-1$, must belong to $QC$.

Consider now factorizations of the form (\ref{2.2}). Suppose that
$\{s_j\}_{0\le j\le m-1}$ is a nonincreasing sequence of nonnegative numbers,  
the matrix functions $W^{\text t}_0, \breve{W}^{\text t}_j, V_0, \breve{V}_j$ 
(see (\ref{2.3})) are thematic, the $u_j$ are unimodular functions such that
the Toeplitz operators $T_{u_j}$ are Fredholm and $\ind T_{u_j}>0$.
Such factorizations are called {\it thematic factorizations}.

It was shown in [PT2] that if $\Phi\in L^\be(\mm_{m,n})$ and 
$\|H_\Phi\|_{\text e}$ is less than the smallest nonzero superoptimal singular
value of $\Phi$, then $\Phi$ is very badly approximable if and only if
it admits a thematic factorization.

The indices $k_j$ of the thematic factorization (\ref{2.2}) ({\it thematic indices}) 
are defined in case $t_j\neq0$: $k_j\df\ind T_{u_j}$ (recall that if $u_j\in QC$,
then $k_j=-\wind u_j$).

It follows from the results of [PY1] that if $\Phi\in L^\be(\mm_{m,n})$
admits a representation (\ref{2.1}) in which $s>0$, $V$ and $W^{\text t}$ are 
thematic matrix functions, $u$ is a unimodular function such that $T_u$ is Fredholm 
with $\ind T_u>0$, and $\|\Psi\|_{L^\be}\le s$, then $\Phi$ is a badly approximable 
matrix function. If $\Phi$ admits a thematic factorization
(\ref{2.2}), then $\Phi$ is very badly approximable with superoptimal
singular values $s_j$, $0\le j\le m-1$ (see [PY1]).

It also follows from the results of [PT2] that if
$\|H_\Phi\|_{\text e}<\|H_\Phi\|$, $r\le\min\{m,n\}$ is such that
$t_{r-1}>\|H_\Phi\|_{\text e}$ and $t_{r-1}>t_r$, and
$F\in \O_{r-1}$, then $\Phi-F$ admits a factorization
\beq
\label{2.4}
\Phi-F=
W_0^*\cdots W^*_{r-1} 
\left( \begin{array}{ccccccc}
t_0 u_0 & 0 & \cdots & 0 & 0 \\
0 & t_1 u_1 & \cdots & 0 & 0 \\
\vdots & \vdots & \ddots & \vdots & \vdots\\
0 & 0 & \cdots & t_{r-1}u_{r-1} & 0\\
0&0&\cdots&0&\Psi\end{array}
\right)
V^*_{r-1}\cdots V^*_0,
\end{equation}
in which the $V_j$ and $W_j$ have the form (\ref{2.3}), the
$W^{\text t}_0, \breve{W}^{\text t}_j, V_0, \breve{V}_j$ are thematic 
matrix functions, the $u_j$ are unimodular functions such that $T_{u_j}$ is
Fredholm and \linebreak$\ind T_{u_j}>0$, 
\beq
\label{2.5}
\|\Psi\|_{L^\be}\le t_{r-1}\quad\mbox{and}\quad\|H_\Psi\|<t_{r-1}.
\end{equation}

Factorizations of the form (\ref{2.4}) with a nonincreasing sequence
$\{t_j\}_{0\le j\le r-1}$ and $\Psi$ satisfying \rf{2.5}
are called {\it partial thematic factorizations}.
If $\Phi-F$ admits a partial thematic factorization of the form (\ref{2.4}),
then $t_0,t_1,\cdots,t_{r-1}$ are the largest $r$ superoptimal singular values of
$\Phi$, and so they do not depend on the choice of a partial thematic factorization.

The matrix entry $\Psi$ in the partial thematic factorization (\ref{2.4}) is called
the {\it residual entry} of the partial thematic factorization.

\

\pagebreak

\setcounter{equation}{0}
\section{\bf Invariance of residual entries}

\

The aim of this section is to show that if a matrix function admits a partial
thematic factorization of the form (\ref{2.4}), then the residual entry $\Psi$ in 
(\ref{2.4}) is uniquely determined by the function itself modulo constant
unitary factors.

\begin{lem}
\label{t3.1}
Let $\Phi$ be an $m\times n$ matrix of the form
$$
\Phi=W^*\left(\begin{array}{cc}u&0\\0&\Psi\end{array}\right)V^*,
$$
where $m,n\ge2$, $u\in\C$, $\Psi\in\mm_{m-1,n-1}$, and 
$$
V=\left(\begin{array}{cc}\bs{v}&\ov{\Theta}\end{array}\right)\in\mm_{n,n},\quad
W=\left(\begin{array}{cc}\bs{w}&\ov{\Xi}\end{array}\right)^{\text t}\in\mm_{m,m}
$$
are unitary matrices such that $\bs{v}\in\mm_{n,1}$ and $\bs{w}\in\mm_{m,1}$. Then
$$
\Psi=\Xi^*\Phi\ov{\Theta}.
$$
\end{lem}

\Pf We have
\begin{eqnarray*}
\Xi^*\Phi\ov{\Theta}&=&
\Xi^*W^*\left(\begin{array}{cc}u&0\\0&\Psi\end{array}\right)V^*\ov{\Theta}\\
&=&\Xi^*\left(\begin{array}{cc}\ov{\bs{w}}&\Xi\end{array}\right)
\left(\begin{array}{cc}u&0\\0&\Psi\end{array}\right)
\left(\begin{array}{c}\bs{v}^*\\\Theta^{\text t}\end{array}\right)\ov{\Theta}\\
&=&\left(\begin{array}{cc}0&I_{m-1}\end{array}\right)
\left(\begin{array}{cc}u&0\\0&\Psi\end{array}\right)
\left(\begin{array}{c}0\\I_{n-1}\end{array}\right)=\Psi.\quad\bl
\end{eqnarray*}

\begin{cor}
\label{t3.2}
Let $\Phi$ be an $m\times n$ matrix of the form
$$
\Phi=W_0^*\cdots W_{r-1}^*
\left(\begin{array}{ccccc}\f_0&0&\cdots&0&0\\
0&\f_1&\cdots&0&0\\
\vdots&\vdots&\ddots&\vdots&\vdots\\
0&0&\cdots&\f_{r-1}&0\\
0&0&\cdots&0&\Psi
\end{array}\right)
V^*_{r-1}\cdots V_0^*,
$$
where $r<\min\{m,n\}$, $~\f_0,\f_1,\cdots,\f_{r-1}\in\C$,
$$
V_j=\left(\begin{array}{cc}I_j&0\\0&\breve V_j\end{array}\right),\quad
W_j=\left(\begin{array}{cc}I_j&0\\0&\breve W_j\end{array}\right),
$$
are unitary matrices such that
{\em$$
\breve V_j=\left(\begin{array}{cc}\bs{v}_j&\ov{\Theta_j}\end{array}\right),\quad
\breve W_j=\left(\begin{array}{cc}\bs{w}_j&\ov{\Xi_j}\end{array}\right)^{\text t},
\quad0\le j\le r-1,
$$}
$\bs{v}_j\in\mm_{n-j,1}$, $\bs{w}_j\in\mm_{m-j,1}$. Then
$$
\Psi=\Xi^*_{r-1}\cdots\Xi^*_1\Xi^*_0\Phi\ov{\Theta_0\Theta_1\cdots\Theta_{r-1}}.
$$
\end{cor}

\Pf The result follows immediately from Lemma \ref{t3.1} by induction. $\bl$

\begin{thm}
\label{t3.3}
Suppose that a matrix function $\Phi\in L^\be(\mm_{m,n})$ admits partial
thematic factorizations
$$
\Phi=W_0^*\cdots W_{r-1}^*
\left(\begin{array}{ccccc}t_0u_0&0&\cdots&0&0\\
0&t_1u_1&\cdots&0&0\\
\vdots&\vdots&\ddots&\vdots&\vdots\\
0&0&\cdots&t_{r-1}u_{r-1}&0\\
0&0&\cdots&0&\Psi
\end{array}\right)
V^*_{r-1}\cdots V_0^*,
$$
and
$$
\Phi=(W_0^\heartsuit)^*\cdots (W_{r-1}^\heartsuit)^*
\left(\begin{array}{ccccc}t_0u^\heartsuit_0&0&\cdots&0&0\\
0&t_1u^\heartsuit_1&\cdots&0&0\\
\vdots&\vdots&\ddots&\vdots&\vdots\\
0&0&\cdots&t_{r-1}u^\heartsuit_{r-1}&0\\
0&0&\cdots&0&\Psi^\heartsuit
\end{array}\right)
(V^\heartsuit_{r-1})^*\cdots(V^\heartsuit_0)^*.
$$
Then there exist constant unitary matrices $U_1\in\mm_{n-r,n-r}$ and 
$U_2\in\mm_{m-r,m-r}$ such that
$$
\Psi^\heartsuit=U_2\Psi U_1.
$$
\end{thm}

Recall that by the definition of a partial thematic factorization, 
$\Psi$ must satisfy \rf{2.5}, and this is very important.

\Pf Let 
$$
V_j=\left(\begin{array}{cc}I_j&0\\0&\breve V_j\end{array}\right),\quad
W_j=\left(\begin{array}{cc}I_j&0\\0&\breve W_j\end{array}\right),
$$
and
$$
V^\heartsuit_j=
\left(\begin{array}{cc}I_j&0\\0&\breve V^\heartsuit_j\end{array}\right),\quad
W^\heartsuit_j=
\left(\begin{array}{cc}I_j&0\\0&\breve W^\heartsuit_j\end{array}\right),
$$
where
$$
\breve V_j=\left(\begin{array}{cc}\bs{v}_j&\ov{\Theta_j}\end{array}\right),\quad
\breve W_j=\left(\begin{array}{cc}\bs{w}_j&\ov{\Xi_j}\end{array}\right)^{\text t},
\quad0\le j\le r-1,
$$
$$
\breve V^\heartsuit_j=
\left(\begin{array}{cc}\bs{v}_j^\heartsuit&\ov{\Theta^\heartsuit_j}\end{array}\right),\quad
\breve W^\heartsuit_j=
\left(\begin{array}{cc}\bs{w}_j^\heartsuit&\ov{\Xi^\heartsuit_j}
\end{array}\right)^{\text t},
\quad0\le j\le r-1.
$$
Here $\breve V_0\df V_0$, $\breve W_0\df V_0$, 
$\breve V_0^\heartsuit\df V_0^\heartsuit$, and
$\breve W_0^\heartsuit\df W_0^\heartsuit$.

We need the following lemma.

\begin{lem}
\label{t3.4}
\beq
\label{3.1}
\Theta_0\Theta_1\cdots\Theta_{r-1}H^2(\C^{n-r})=
\Theta^\heartsuit_0\Theta^\heartsuit_1\cdots\Theta^\heartsuit_{r-1}H^2(\C^{n-r})
\end{equation}
and
\beq
\label{3.2}
\Xi_0\Xi_1\cdots\Xi_{r-1}H^2(\C^{m-r})=
\Xi^\heartsuit_0\Xi^\heartsuit_1\cdots\Xi^\heartsuit_{r-1}H^2(\C^{m-r}).
\end{equation}
\end{lem}

Let us first complete the proof of Theorem \ref{t3.3}. Consider the inner 
matrix functions 
$$
\Theta=\Theta_0\Theta_1\cdots\Theta_{r-1},\quad
\Theta^\heartsuit=
\Theta^\heartsuit_0\Theta^\heartsuit_1\cdots\Theta^\heartsuit_{r-1}
$$
and
$$
\Xi=\Xi_0\Xi_1\cdots\Xi_{r-1},\quad
\Xi^\heartsuit=
\Xi^\heartsuit_0\Xi^\heartsuit_1\cdots\Xi^\heartsuit_{r-1}.
$$
By Lemma \ref{t3.4}, $\Theta H^2(\C^{n-r})=\Theta^\heartsuit H^2(\C^{n-r})$.
It is well known that in this case there exists a constant unitary matrix
$Q_1\in\mm_{n-r,n-r}$ such that $\Theta^\heartsuit=\Theta Q_1$ ($\Theta$ and
$\Theta^\heartsuit$ determine the same invariant subspace under multiplication
by $z$, see e.g., [N]). Similarly, there exists a constant unitary matrix
$Q_2\in\mm_{m-r,m-r}$ such that $\Xi^\heartsuit=\Xi Q_2$.

By Corollary \ref{t3.2},
$$
\Psi=\Xi^*\Phi\ov{\Theta},\quad
\Psi^\heartsuit=(\Xi^\heartsuit)^*\Phi\ov{\Theta^\heartsuit}.
$$
Hence,
$$
\Psi^\heartsuit=Q_2^*\Xi^*\Phi\ov{\Theta}\ov{Q_1}=
Q_2^*\Psi\ov{Q_1}.\quad\bl
$$

{\bf Proof of Lemma \ref{t3.4}.} It is sufficient to prove (\ref{3.1}). Indeed,
(\ref{3.2}) follows from (\ref{3.1}) applied to $\Phi^{\text t}$.

It is easy to see that without loss of generality we may assume that \linebreak
$\|\Psi\|_{L^\be}<t_{r-1}$. Indeed, we can subtract from $\Phi$ a matrix function in
$\O_{r-1}$, and it follows from Lemma 1.5 of [PY1] that the resulting function
admits a partial thematic factorization with the same unitary-valued function
$V_j$ and $W_j$, $0\le j\le r-1$, and residual entry whose $L^\be$ norm is less
that $t_{r-1}$. It is also easy to see that if $\|\Psi\|_{L^\be}<t_{r-1}$, then
$\|\Psi^\heartsuit\|_{L^\be}$ must also be less than $t_{r-1}$.

Consider the subspace $\cL$ of $H^2(\C^n)$ defined by
$$
\cL=\left\{f\in H^2(\C^n):~V^{\text t}_{r-1}\cdots V^{\text t}_1V^{\text t}_0f=
\left(
\begin{array}{c}0\\\vdots\\0\\\ast\\\vdots\\\ast\end{array}\right)
\!\!\!\!\!\!\!\!\!
\begin{array}{c}
\left.\begin{array}{c}\\\\\\\end{array}\right\}r
\\\\\\\\
\end{array}\right\},
$$
i.e., $\cL$ consists of vector functions $f\in H^2(\C^n)$ such that the first
$r$ components of the vector function 
$V^{\text t}_{r-1}\cdots V^{\text t}_1V^{\text t}_0f$ are zero.

We define the real function $\r$ on $\R$ by
$$
\r(x)=\left\{\begin{array}{ll}x,&x\ge t^2_{r-1}\\0,&x<t^2_{r-1}\end{array}
\right.
$$
and consider the operator $M:H^2(\C^n)\to L^2(\C^n)$ of multiplication by 
the matrix function $\r(\Phi^{\text t}\ov{\Phi})$:
$$
Mf=\r(\Phi^{\text t}\ov{\Phi})f,\quad f\in H^2(\C^n).
$$
Let us show that
\beq
\label{3.3}
\cL=\Ker M.
\end{equation}

We have
$$
\Phi^{\text t}\ov{\Phi}=\ov{V_0V_1\cdots V_{r-1}}
\left(\begin{array}{ccccc}t_0^2&\cdots&0&0\\
\vdots&\ddots&\vdots&\vdots\\
0&\cdots&t_{r-1}^2&0\\
0&\cdots&0&\Psi^{\text t}\ov{\Psi}\end{array}\right)
V^{\text t}_{r-1}\cdots V^{\text t}_1V^{\text t}_0,
$$
and since $\|\Psi^{\text t}\ov{\Psi}\|_{L^\be}<t_{r-1}^2$, it follows that
$$
\r(\Phi^{\text t}\ov{\Phi})=\ov{V_0V_1\cdots V_{r-1}}
\left(\begin{array}{ccccc}t_0^2&\cdots&0&0\\
\vdots&\ddots&\vdots&\vdots\\
0&\cdots&t_{r-1}^2&0\\
0&\cdots&0&0\end{array}\right)
V^{\text t}_{r-1}\cdots V^{\text t}_1V^{\text t}_0.
$$
Since all matrix functions $V_j$ are unitary-valued, this implies (\ref{3.3}).

Thus the subspace  $\cL$ is uniquely determined by the function $\Phi$ and
does not depend on the choice of a partial thematic factorization. It is easy 
to see that to complete the proof of Lemma \ref{t3.4}, it is sufficient
to prove the following lemma.

\begin{lem}
\label{t3.5}
\beq
\label{3.4}
\cL=\Theta_0\Theta_1\cdots\Theta_{r-1}H^2(\C^{n-r}).
\end{equation}
\end{lem}

\Pf We show by induction on $r$ that (\ref{3.4}) holds even without the assumption
that $\|\Psi\|_{L^\be}<t_{r-1}$ (note that this assumption is very important in the
proof of (\ref{3.3})). 

Suppose that $r=1$. Then
$$
\cL=\left\{f\in H^2(\C^n):~V_0^{\text t}f=
\left(\begin{array}{c}0\\\ast\\\vdots\\\ast\end{array}\right)\right\}.
$$
Obviously, if $f\in\Theta_0H^2(\C^{n-1})$, then $f\in\cL$. Suppose now that
$f\in\cL$. We have
$$
V_0^{\text t}f=\left(\begin{array}{c}0\\g\end{array}\right),
\quad g\in L^2(\C^{n-1}).
$$
Then
$$
f=\ov{V_0}\left(\begin{array}{c}0\\g\end{array}\right)
=\left(\begin{array}{cc}\ov{\bs{v}_0}&\Theta_0\end{array}\right)
\left(\begin{array}{c}0\\g\end{array}\right)=\Theta_0g.
$$
Let us show that $g\in H^2(\C^{n-1})$. It suffices to prove that
$g^{\text t}\g\in H^2$ for any constant vector $\g\in\C^{n-1}$.
Since $\Theta_0^{\text t}$ is outer, there exists a sequence $\{\f_n\}_{n\ge0}$
of functions in $H^2(\C^n)$ such that
$$
\lim_{n\to\be}\Theta_0^{\text t}\f_n\to\g\quad\mbox{in}\quad H^2(\C^{n-1}).
$$
We have
$$
f^{\text t}\f_n=g^{\text t}\Theta_0^{\text t}\f_n\to g^{\text t}\g
\quad\mbox{in}\quad H^1,
$$
and so $g^{\text t}\g\in H^2$ which proves the result for $r=1$.

Suppose now that $r\ge2$. By the induction hypothesis
$$
\cL=\left\{\Theta_0\cdots\Theta_{r-2}g:~g\in H^2(\C^{n-r+1}),~
V^{\text t}_{r-1}\cdots V^{\text t}_0\Theta_0\cdots\Theta_{r-2}g=
\left(
\begin{array}{c}0\\\vdots\\0\\\ast\\\vdots\\\ast\end{array}\right)
\!\!\!\!\!\!\!\!\!\!
\begin{array}{c}
\left.\begin{array}{c}\\\\\\\end{array}\right\}r
\\\\\\\\
\end{array}
\right\}.
$$
It follows from the definition of thematic matrix functions that
$$
V^{\text t}_{r-2}\cdots V^{\text t}_0\Theta_0\cdots\Theta_{r-2}
=\left(\begin{array}{cccc}0&\cdots&0\\
\vdots&\ddots&\vdots\\
0&\cdots&0\\
1&\cdots&0\\
\vdots&\ddots&\vdots\\
0&\cdots&1
\end{array}\right)
\!\!\!\!\!\!\!
\begin{array}{c}
\!\!\!\!\!\!\!\!\!\!\!\left.\begin{array}{c}\\\\\\\end{array}\right\}r-1\\
\left.\begin{array}{c}\\\\\\\end{array}\right\}n-r+1
\end{array}
=\left(\begin{array}{c}0\\I_{n-r+1}\end{array}\right).
$$
Hence,
$$
\cL=\left\{\Theta_0\cdots\Theta_{r-2}g:~g\in H^2(\C^{n-r+1}),~
\left(\begin{array}{c}\bs{v}_{r-1}^{\text t}\\\Theta^*_{r-1}\end{array}\right)g=
\left(\begin{array}{c}0\\\ast\\\vdots\\\ast\end{array}\right)
\right\}.
$$
Since the result has already been proved for $r=1$,
$$
\cL=\left\{\Theta_0\cdots\Theta_{r-2}g:~g\in\Theta_{r-1}H^2(\C^{n-r})\right\}
=\Theta_0\cdots\Theta_{r-1}H^2(\C^{n-r}).\quad\bl
$$

\

\setcounter{equation}{0}
\section{\bf Monotone thematic factorizations and invariance of indices}

\

In this section we study the problem of the invariance of indices of thematic 
factorizations of very badly approximable matrix functions. In [PY1] it was shown 
that the indices of a thematic factorization are not determined uniquely by the 
matrix function but may depend on the choice of a thematic factorization. 
For example, the matrix function 
$\Phi=\left(\begin{array}{cc}{\bar z}^2&0\\0&\bar z^6\end{array}\right)$
admits the following thematic factorizations

\begin{eqnarray*}
\Phi&=&\left(\begin{array}{cc}1&0\\0&1\end{array}\right)
\left(\begin{array}{cc}{\bar z}^2&0\\0&\bar z^6\end{array}\right)
\left(\begin{array}{cc}1&0\\0&1\end{array}\right)\\[.5pc]
&=&\left(\begin{array}{cc}\frac{1}{\sqrt{2}}&-\frac{z^5}{\sqrt{2}}\\[.8pc]
\frac{\bar z^5}{\sqrt{2}}&\frac{1}{\sqrt{2}}\end{array}\right)
\left(\begin{array}{cc}{\bar z}&0\\[.8pc]0&\bar z^7\end{array}\right)
\left(\begin{array}{cc}\frac{\bar z}{\sqrt{2}}&\frac{1}{\sqrt{2}}\\[.8pc]
-\frac{1}{\sqrt{2}}&\frac{z}{\sqrt{2}}\end{array}\right)\\[.5pc]
&=&\left(\begin{array}{cc}0&1\\1&0\end{array}\right)
\left(\begin{array}{cc}{\bar z}^6&0\\0&\bar z^2\end{array}\right)
\left(\begin{array}{cc}0&1\\1&0\end{array}\right).
\end{eqnarray*}
The superoptimal singular values of $\Phi$ are $t_0=t_1=1$.
The indices of the first factorization are 2, 6, the indices of the second are 1, 7, 
and the indices of the third are 6, 2. Note that
for all above factorizations the sum of the indices is 8.

In [PY2] it was shown (in the case of $H^\be+C$ functions) that the sum of thematic 
indices that correspond to all superoptimal singular values equal to a positive 
specific value does not depend on the choice of a thematic factorization. In other 
words, for each positive superoptimal singular value $t$ the numbers
$$
\nu_t\df\sum_{\{j:t_j=t\}}k_j
$$
do not depend on the choice of a thematic factorization. The same result was obtained
in [PT2] in the case when $\|H_\Phi\|_{\text e}$ is less than the smallest
nonzero superoptimal singular value. Note that it also follows from
the results of [PY2] and [PT2] that the same invariance property holds for
partial thematic factorizations.

A natural question arises of whether we can distribute arbitrarily the numbers
$\nu_t$ between the indices $k_j$ with $t_j=t$  by choosing an appropriate thematic
factorization (recall that the $k_j$ must be positive integers).

In this section we show that the answer to this question is negative. 

{\bf Definition.} A (partial) thematic factorization is called {\it monotone} 
if for any positive superoptimal singular value $t$ the thematic indices 
$k_r,k_{r+1},\cdots,k_s$ that correspond to all superoptimal singular values 
equal to $t$ satisfy
\beq
\label{4.1}
k_r\ge k_{r+1}\ge\cdots\ge k_s.
\end{equation}
Here $t_r,t_{r+1},\cdots,t_s$ are the superoptimal singular values equal to $t$.

We prove in this section that if $\|H_\Phi\|_{\text e}$ is less than the
smallest nonzero superoptimal singular value of $\Phi$, then $\Phi-\A\Phi$ possesses
a monotone thematic factorization. We also show that the indices
of a monotone thematic factorization are uniquely determined by the function 
$\Phi$ itself and do not depend on the choice of a thematic factorization. In
particular this is the case if $\Phi\in(H^\be+C)(\mm_{m,n})$. The same results
also hold for partial thematic factorizations.

In the above example only the third thematic factorization is monotone. It will 
follow from the results of this section that the thematic indices of any monotone
thematic factorization must be equal to 6, 2. In particular, there are no thematic 
factorizations with indices 7, 1. Note that it is important that the
indices in (\ref{4.1}) are arranged in the {\it nonincreasing} order. The above 
example shows that the first two thematic factorizations have different
thematic indices 2, 6 and 1, 7 that are arranged in the increasing order.

\begin{thm}
\label{t4.1}
Suppose that $\Phi\in L^\be(\mm_{m,n})$ and $r\le\min\{m,n\}$ is a positive integer
such that the superoptimal singular values of $\Phi$ satisfy {\em
$$
t_{r-1}>t_r,\quad t_{r-1}>\|H_\Phi\|_{\text e}.
$$}
If $\Phi$ admits a partial thematic factorization of the form {\em(\ref{2.4})},
then $\Phi$ admits a monotone partial thematic factorization of the form 
{\em(\ref{2.4})}.
\end{thm}

\Pf Clearly, $\|H_{z^j\Phi}\|=\dist_{L^\be}(\Phi,\bar z^jH^\be(\mm_{m,n}))$, 
and it is easy to see that
$$
\lim_{j\to\be}\|H_{z^j\Phi}\|=\dist_{L^\be}\big(\Phi,(H^\be+C)(\mm_{m,n})\big)=
\|H_\Phi\|_{\text e}<\|H_\Phi\|.
$$

Put
$$
\iota(H_\Phi)\df\min\{j\ge0:~\|H_{z^j\Phi}\|<\|H_\Phi\|\}.
$$
Obviously, $\iota(H_\Phi)$ depends only on the Hankel operator $H_\Phi$ and does not
depend on the choice of its symbol. 

We need three lemmas.

\begin{lem}
\label{t4.2}
Let $\Phi$ be a matrix function in $L^\be(\mm_{m,n})$ such that
{\em$\|H_\Phi\|_{\text e}<\|H_\Phi\|$}. Suppose that
\beq
\label{4.2}
\Phi=W^*\left(\begin{array}{cc}tu&0\\0&\U\end{array}\right)V^*,
\end{equation}
where $V$ and {\em$W^{\text t}$} are thematic matrix functions of sizes 
$n\times n$ and $m\times m$, $t>0$, $\|\U\|_{L^\be}\le t$, 
and $u$ is a unimodular function such that $T_u$ is Fredholm.
Then $\ind T_u\le\iota(H_\Phi)$.
\end{lem}

\begin{lem}
\label{t4.3}
Let $\Phi$ be a badly approximable matrix function in 
$L^\be(\mm_{m,n})$ such that {\em$\|H_\Phi\|_{\text e}<\|H_\Phi\|$}. 
Then $\Phi$ admits a representation {\em(\ref{4.2})}
with thematic matrix functions $V$ and {\em$W^{\text t}$},
$t=t_0=\|H_\Phi\|$, and a unimodular function $u$ such that $T_u$ is Fredholm and
$$
\ind T_u=\iota(H_\Phi).
$$
\end{lem}

\begin{lem}
\label{t4.4}
Let $\Phi\in L^\be(\mm_{m,n})$ be a matrix function of the form
$$
\Phi=W^*\left(\begin{array}{cc}u&0\\0&\U\end{array}\right)V^*,
$$
where $V$ and {\em$W^{\text t}$} are thematic matrix functions of sizes 
$n\times n$ and $m\times m$, $u$ is a unimodular function such that $T_u$ is
Fredholm, $\ind T_u=0$, $\|H_\U\|\le1$, and {\em$\|H_\U\|_{\text e}<1$}. 
If $\|H_\Phi\|<1$, then $\|H_\U\|<1$.
\end{lem}

Let us first complete the proof of Theorem \ref{t4.1}. We argue by induction on $r$. 
For $r=1$ the result is trivial. Suppose now that $r>1$. By Lemma \ref{t4.3},
$\Phi$ admits a representation
$$
\Phi=W^*\left(\begin{array}{cc}t_0u_0&0\\0&\U\end{array}\right)V^*,
$$
where $V$ and $W^{\text t}$ are thematic functions, $\|\U\|_{L^\be}\le t_0$, 
and $u_0$ is a unimodular
function such that $T_{u_0}$ is Fredholm and $\ind T_{u_0}=\iota(H_\Phi)$. By
Theorem 6.3 of [PT2], 
\beq
\label{4.3}
\|H_\U\|_{\text e}\le\|H_\Phi\|_{\text e}.
\end{equation}
It follows from the results of \S 4 and \S 6 of [PT2] that $\U$ admits a partial 
thematic factorization of the form
$$
\U=W_1^*\cdots W_{r-1}^*\left(\begin{array}{ccccc}
t_1u_1&\cdots&0&0\\
\vdots&\ddots&\vdots&\vdots\\
0&\cdots&t_{r-1}u_{r-1}&0\\
0&\cdots&0&\Psi
\end{array}\right)V_{r-1}^*\cdots V_1^*.
$$
By the induction hypothesis we may assume that this partial thematic factorization 
is monotone. Clearly, $t_1=\|\U\|_{L^\be}$. 
If $t_1<t_0$, then it is easy to see that the above factorization of $\U$ leads
to a monotone partial thematic factorization of $\Phi$.

Suppose now that $t_1=t_0$.  To prove that the above factorization of $\U$
leads to a monotone partial thematic factorization of $\Phi$,
we have to establish the inequality 
$\ind T_{u_0}\ge\ind T_{u_1}$. By Lemma \ref{t4.2}, 
$\iota(H_\U)\ge\ind T_{u_1}$, and it suffices to prove the inequality 
$$
\iota(H_\Phi)=\ind T_{u_0}\ge\iota(H_\U).
$$
Put $\iota\df\iota(H_\Phi)$. We have
$$
z^\iota\Phi=W^*\left(\begin{array}{cc}t_0z^\iota u_0&0\\0&z^\iota\U
\end{array}\right)V^*.
$$
Clearly, $\ind T_{z^\iota u_0}=0$.
By the definition of $\iota$, $\|H_{z^\iota\Phi}\|<\|H_\Phi\|=t_0$. 
It is easy to see that 
$$
\|H_{z^\iota\U}\|_{\text e}=\|H_{\U}\|_{\text e}<t_0 
$$
by (\ref{4.3}). It follows from Lemma \ref{t4.4} that $\|H_{z^\iota\U}\|<t_0$ 
which means that $\iota(H_\U)\le\iota$. $\bl$

{\bf Proof of Lemma \ref{t4.2}.} Let $k=\ind T_u$. Clearly, it is sufficient to 
consider the case $k>0$. Then $\Phi$ is badly approximable and $\|H_\Phi\|=t$
(see \S 2). We have
$$
z^{k-1}\Phi=W^*\left(\begin{array}{cc}tz^{k-1}u&0\\0&z^{k-1}\U\end{array}\right)V^*.
$$
Then $\wind(z^{k-1}u)=-1$, and so $z^{k-1}\Phi$ is badly approximable 
and $\|\Phi\|_{L^\be}=t$ (see \S 2). Hence, 
$$
\|H_{z^{k-1}\Phi}\|=\|z^{k-1}\Phi\|_{L^\be}=\|\Phi\|_{L^\be}=t=\|H_\Phi\|,
$$
and so $\iota(H_\Phi)\ge k$. $\bl$

{\bf Proof of Lemma \ref{t4.3}.} Put $\iota\df\iota(H_\Phi)$. Then 
$$
\|H_{z^{\iota-1}\Phi}\|=\|H_\Phi\|=\|\Phi\|_{L^\be}=\|z^{\iota-1}\Phi\|_{L^\be},
$$
and so $z^{\iota-1}\Phi$ is badly approximable. Clearly,
$$
\|H_{z^{\iota-1}\Phi}\|_{\text e}=\|H_\Phi\|_{\text e}<
\|H_\Phi\|=\|H_{z^{\iota-1}\Phi}\|.
$$
Hence, (see \S 2) $z^{\iota-1}\Phi$ admits a representation
$$
z^{\iota-1}\Phi=W^*\left(\begin{array}{cc}t\o&0\\0&\O\end{array}\right)V^*,
$$
where $t=\|H_\Phi\|$, $\o$ is a unimodular function such that $\ind T_\o>0$,
$V$ and $W^{\text t}$ are thematic functions and $\|\O\|_{L^\be}\le t$. 
Therefore
$$
\Phi=W^*\left(\begin{array}{cc}t\bar z^{\iota-1}\o&0\\0&\bar z^{\iota-1}\O
\end{array}\right)V^*.
$$
Let $u=\bar z^{\iota-1}\o$. Clearly, $\ind T_u\ge\iota$. Finally, by Lemma \ref{t4.2},
$\ind T_u=\iota$. $\bl$

{\bf Proof of Lemma \ref{t4.4}.} The proof is based on the argument given
in the proof of Lemma 1.2 of [PY2]. Let
$$
V=\left(\begin{array}{cc}\bs{v}&\ov{\Theta}\end{array}\right),\quad
W^{\text t}=\left(\begin{array}{cc}\bs{w}&\ov{\Xi}\end{array}\right).
$$
By Theorem 5.1 of [PT2], there exist $A\in H^\be(\mm_{n-1,n})$ and 
$B\in H^\be(\mm_{m-1,m})$ such that $A\Theta=I_{n-1}$ and $B\Xi=I_{m-1}$.
Without loss of generality we may assume that $\|\U\|_{L^\be}\le1$.

Suppose that $\|H_\U\|=1$. Since $\|H_\U\|_{\text e}<1$, there exists
a nonzero function \linebreak$g\in H^2(\C^{n-1})$ such that 
$\|H_\U g\|_2=\|g\|_2$. Then $\U g\in H^2_-(\C^{m-1})$ and \linebreak
$\|\U(\z)g(\z)\|_{\C^{m-1}}=\|g(\z)\|_{\C^{n-1}}$ for almost all $\z\in\T$.

Let 
$$
f=A^{\text t}g+\bs{v}q,
$$
where $q$ is a scalar function in $H^2$. We want to find such a $q$ that
$\|H_\Phi f\|_2=\|f\|_2$. Note that $f$ is a nonzero function
since
$$
V^*f=\left(\begin{array}{c}\bs{v}^*\\\Theta^{\text t}\end{array}\right)
(A^{\text t}g+\bs{v}q)=
\left(\begin{array}{c}\bs{v}^*A^{\text t}g+q\\g\end{array}\right)
$$
and $g\neq0$.

We have
\begin{eqnarray*}
\Phi f&=&W^*\left(\begin{array}{cc}u&0\\0&\U\end{array}\right)
\left(\begin{array}{c}\bs{v}^*A^{\text t}g+q\\g\end{array}\right)\\
&=&\left(\begin{array}{cc}\ov{\bs{w}}&\Xi\end{array}\right)
\left(\begin{array}{c}u\bs{v}^*A^{\text t}g+uq\\\U g\end{array}\right)\\
&=&\ov{\bs{w}}(u\bs{v}^*A^{\text t}g+uq)+\Xi\U g.
\end{eqnarray*}
Since the matrix functions $W^*$ and $V^*$ are unitary-valued and 
$\|\U(\z)g(\z)\|_{\C^{m-1}}=\|g(\z)\|_{\C^{n-1}}$, it follows that
$\|\Phi(\z)f(\z)\|_{\C^m}=\|f(\z)\|_{\C^n}$. It remains to choose $q$ so that
$\Phi f\in H^2_-(\C^m)$.

Since $W^*$ is a unitary-valued matrix function, we have
$$
I_m=\left(\begin{array}{cc}\ov{\bs{w}}&\Xi\end{array}\right)
\left(\begin{array}{c}\bs{w}^{\text t}\\\Xi^*\end{array}\right)=
\ov{\bs{w}}\bs{w}^{\text t}+\Xi\Xi^*.
$$
Hence, 
$$
\Xi=\Xi(B\Xi)^*=\Xi\Xi^*B^*=(I_m-\ov{\bs{w}}\bs{w}^{\text t})B^*.
$$
It follows that
\begin{eqnarray*}
\Phi f&=&
\ov{\bs{w}}(u\bs{v}^*A^{\text t}g+uq)+(I_m-\ov{\bs{w}}\bs{w}^{\text t})B^*\U g\\
&=&\ov{\bs{w}}(u\bs{v}^*A^{\text t}g+uq-\bs{w}^{\text t}B^*\U g)+B^*\U g.
\end{eqnarray*}
Clearly, $B^*\U g\in H^2_-(\C^m)$, and so it suffices to find $q\in H^2$
such that
$$
u\bs{v}^*A^{\text t}g+uq-\bs{w}^{\text t}B^*\U g\in H^2_-
$$
which is equivalent to the condition
$$
T_uq=\pp_+(\bs{w}^{\text t}B^*\U g-u\bs{v}^*A^{\text t}g).
$$
The existence of such a $q$ follows from the well-known fact that the Toeplitz 
operator $T_u$ is invertible; indeed, it is Fredholm and $\ind T_u=0$
(see e.g., [D] or [N]). $\bl$

\begin{cor}
\label{t4.5}
Let $\Phi$ be a very badly approximable matrix function in
$L^\be(\mm_{m,n})$ such that {\em$\|H_\Phi\|_{\text e}$} is less than
the smallest nonzero superoptimal singular value of $\Phi$. 
Then $\Phi$ admits a monotone thematic factorization.
\end{cor}

\begin{cor}
\label{t4.6}
Let $\Phi$ be a very badly approximable matrix function in \linebreak
$(H^\be+C)(\mm_{m,n})$. Then $\Phi$ admits a monotone thematic factorization.
\end{cor}

We are going to prove now that the indices of a monotone thematic factorization
are uniquely determined by the function itself. We need the following lemma.

\begin{lem}
\label{t4.7}
Suppose that a matrix function $\Phi\in L^\be(\mm_{m,n})$ admits a factorization of 
the form
$$
\Phi=W^*_0\cdots W^*_{r-1}
\left(\begin{array}{ccccc}tu_0&0&\cdots&0&0\\0&tu_1&\cdots&0&0\\
\vdots&\vdots&\ddots&\vdots&\vdots\\0&0&\cdots&tu_{r-1}&0\\
0&0&\cdots&0&\Psi
\end{array}\right)
V^*_{r-1}\cdots V_0^*,
$$
where the $V_j$ and $W_j$ are of the form {\em(\ref{2.3})},
$\|H_\Psi\|<t$ and the $u_j$ are unimodular functions such that $T_{u_j}$ is
Fredholm and $\ind T_{u_j}\le0$. 
If {\em$\|H_\Phi\|_{\text e}<t$}, then $\|H_\Phi\|<t$.
\end{lem}

\Pf We argue by induction on $r$. Let $r=1$.
We have 
$$
\Phi=W^*\left(\begin{array}{cc}tu&0\\0&\Psi\end{array}\right)V^*,
$$
where $V$ and $W^{\text t}$ are thematic matrix functions,
$u$ is a unimodular function such that $T_u$ is Fredholm, $\ind T_u\le0$, 
and $\|H_\Psi\|<t$. It follows from Lemma 1.5 of [PY1] that we may subtract from 
$\Psi$ a best analytic approximation without changing $H_\Phi$, and so we may assume that 
$\|\Psi\|_{L^\be}<t$. Without loss of generality we may also assume that $t=1$.

Suppose that $\|H_\Phi\|=1$. Since 
$\|H_\Phi\|_{\text e}<1$, there exists a nonzero 
function $f\in H^2(\C^n)$ such that $\|H_\Phi f\|_2=\|f\|_2$. 
Then $\|\Phi f\|_2=\|f\|_2$ and since $\|\Psi\|_{L^\be}<1$, it follows that $V^*f$
has the form
\beq
\label{4.4}
V^*f=\left(\begin{array}{c}*\\0\\\vdots\\0\end{array}\right).
\end{equation}
Let $\bs{v}$ be the first column of $V$. Equality (\ref{4.4}) means that for almost
all $\z\in\T$ the remaining columns of $V(\z)$ are orthogonal to $f(\z)$ in $\C^n$.
Since $V$ is unitary-valued, it follows that $f=\xi \bs{v}$ for a scalar function
$\xi\in L^2$. Using the fact that $\bs{v}$ is co-outer, we can find a sequence
of $n\times1$ functions $\f_j$ in $H^2$ such that 
${\displaystyle\lim_{j\to\be}}\|\f_j^{\text t}\bs{v}-1\|_2=0$. 
Hence, $\xi$ is the limit in $L^1$ of
the sequence $\f^{\text t}_j f$, and so $\xi\in H^2$. Note that 
$\|f\|_{H^2(\C^n)}=\|\xi\|_{H^2}$.

We have 
$$
\Phi f=W^*\left(\begin{array}{c}u\xi\\0\\\vdots\\0\end{array}\right)
=u\xi\ov{\bs{w}},
$$
where $\bs{w}$ is the first column of $W^{\text t}$. Since $f$ is a maximizing vector
of $H_\Phi$, we have $u\xi\ov{\bs{w}}\in H^2_-(\C^n)$. Again, using the fact that
$\bs{w}$ is co-outer, we find that $u\xi\in H^2_-$, i.e., $\xi\in\Ker T_u$. However,
$T_u$ has trivial kernel since $\ind T_u\le0$. We have got a contradiction. 

Suppose now that $r>1$. Again, we may assume that $\|\Psi\|_{L^\be}<t$. 
Let $d$ be a negative integer such that $d<\ind T_{u_j}$, $0\le j\le r-1$.
Then
$$
z^d\Phi=W^*_0\cdots W^*_{r-1}
\left(\begin{array}{ccccc}tz^du_0&0&\cdots&0&0\\0&tz^du_1&\cdots&0&0\\
\vdots&\vdots&\ddots&\vdots&\vdots\\0&0&\cdots&tz^du_{r-1}&0\\
0&0&\cdots&0&z^d\Psi
\end{array}\right)
V^*_{r-1}\cdots V_0^*
$$
is a partial thematic factorization of $z^d\Phi$. Put
$$
\U=W^*_1\cdots W^*_{r-1}
\left(\begin{array}{ccccc}tu_1&0&\cdots&0&0\\
\vdots&\vdots&\ddots&\vdots&\vdots\\0&0&\cdots&tu_{r-1}&0\\
0&0&\cdots&0&\Psi
\end{array}\right)
V^*_{r-1}\cdots V_1^*
$$
Since obviously, $\|H_{z^d\U}\|_{\text e}=\|H_\U\|_{\text e}$ for any 
$d\in\Z$, it follows from
Theorem 6.3 of [PT2] that $\|H_{z^d\U}\|_{\text e}<t$, and so by the induction 
hypotheses, $\|H_\U\|<t$. We have 
$$
\Phi=W_0^*\left(\begin{array}{cc}tu&0\\0&\U\end{array}\right)V_0^*.
$$
The result follows now from the case $r=1$ which has already been established.
$\bl$

\begin{thm}
\label{t4.8}
Let $\Phi$ be a  badly approximable function in $L^\be(\mm_{m,n})$ 
such that {\em$\|H_\Phi\|_{\text e}<\|H_\Phi\|$}
and let $r$ be the number of superoptimal singular values of $\Phi$ equal to 
$t_0=\|H_\Phi\|$. Consider a monotone partial thematic factorization of $\Phi$
with indices
\beq
\label{4.5}
k_0\ge\cdots\ge k_{r-1}
\end{equation}
corresponding to the superoptimal singular values equal to $t_0$. Let $\kappa\ge0$.
Then
\beq
\label{4.6}
\dim\{f\in H^2(\C^n):\|H_{z^\kappa\Phi}f\|_2=t_0\|f\|_2\}=
\sum_{\{j\in[0,r-1]:k_j>\kappa\}}k_j-\kappa.
\end{equation}
\end{thm}

\Pf Let
$$
\Phi=\left(\begin{array}{ccccc}t_0u_0&0&\cdots&0&0\\
0&t_0u_1&\cdots&0&0\\
\vdots&\vdots&\ddots&\vdots&\vdots\\
0&0&\cdots&t_0u_{r-1}&0\\
0&0&\cdots&0&\Psi\end{array}\right)
$$
be a partial thematic factorization of $\Phi$ with indices satisfying (\ref{4.5}). 
If $\kappa\ge k_0$, then (\ref{4.6}) holds by Lemma \ref{t4.7}. Suppose now that 
$\kappa<k_0$. Let
$$
q=\max\{j\in[0,r-1]:~k_j>\kappa\}.
$$
Clearly, the function $z^\kappa\Phi$ admits the following representation
$$
z^\kappa\Phi=W_0^*\cdots W^*_q
\left(\begin{array}{ccccc}
t_0z^\kappa u_0&0&\cdots&0&0\\
0&t_0z^\kappa u_1&\cdots&0&0\\
\vdots&\vdots&\ddots&\vdots&\vdots\\
0&0&\cdots&t_0z^\kappa u_q&0\\
0&0&\cdots&0&\U
\end{array}\right)V^*_q\cdots V_0^*.
$$
where $\U$ is a matrix function satisfying 
the hypotheses of Lemma \ref{t4.7}. By Lemma \ref{t4.7}, $\|H_\U\|<t_0$. Let 
$R\in H^\be$ be a matrix function such that $\|\U-R\|_{L^\be}<t_0$.
It is easy to show by induction on $q$ that if we perturb $\U$ by a bounded
analytic matrix function, $z^\kappa\Phi$ also changes by an analytic matrix function
(this is the trivial part of Lemma 1.5 of [PY1]). In particular, we can find
a matrix function $G\in H^\be$ such that
$$
z^\kappa\Phi-G=W_0^*\cdots W^*_q
\left(\begin{array}{ccccc}
t_0z^\kappa u_0&0&\cdots&0&0\\
0&t_0z^\kappa u_1&\cdots&0&0\\
\vdots&\vdots&\ddots&\vdots&\vdots\\
0&0&\cdots&t_0z^\kappa u_q&0\\
0&0&\cdots&0&\U-R
\end{array}\right)V^*_q\cdots V_0^*.
$$

By Theorem 9.3 of [PT2], 
$$
\dim\{f\in H^2(\C^n):~\|H_{z^\kappa\Phi-G}f\|_2=t_0\|f\|_2\}=
\sum_{\{j\in[0,r-1]:k_j>\kappa\}}k_j-\kappa
$$
(this equality was stated in [PT2] for thematic factorizations but the
same proof also works for partial thematic factorizations).
Equality (\ref{4.6}) follows now from the obvious fact that 
$H_{z^\kappa\Phi-G}=H_{z^\kappa\Phi}$. $\bl$

We can now deduce from (\ref{4.6}) the following result.

\begin{thm}
\label{t4.9}
Suppose that $\Phi\in L^\be(\mm_{m,n})$ and $q\le\min\{m,n\}$ is a positive integer
such that the superoptimal singular values of $\Phi$ satisfy {\em
$$
t_{q-1}>t_q,\quad t_{q-1}>\|H_\Phi\|_{\text e}
$$}
and $\Phi$ admits a monotone partial thematic factorization
$$
\Phi=W_0^*\cdots W_{q-1}^*\left(\begin{array}{ccccc}
t_0u_0&0&\cdots&0&0\\
0&t_1u_1&\cdots&0&0\\
\vdots&\vdots&\ddots&\vdots&\vdots\\
0&0&\cdots&t_{q-1}u_{q-1}&0\\
0&0&\cdots&0&\U
\end{array}\right)V^*_{q-1}\cdots V^*_0.
$$
Then the indices of this factorization  
are uniquely determined by the function $\Phi$ itself.
\end{thm}

\Pf Let $r$ be the number of superoptimal singular values equal to 
$\|H_\Phi\|$. Then $\Phi$ admits the following partial thematic factorization
$$
\Phi=W_0^*\cdots W_{r-1}^*\left(\begin{array}{ccccc}
t_0u_0&0&\cdots&0&0\\
0&t_1u_1&\cdots&0&0\\
\vdots&\vdots&\ddots&\vdots&\vdots\\
0&0&\cdots&t_{r-1}u_{r-1}&0\\
0&0&\cdots&0&\Psi
\end{array}\right)V^*_{r-1}\cdots V^*_0,
$$
where
$$
\Psi=W_r^*\cdots W_{q-1}^*\left(\begin{array}{ccccc}
t_ru_r&\cdots&0&0\\
\vdots&\ddots&\vdots&\vdots\\
0&\cdots&t_{q-1}u_{q-1}&0\\
0&\cdots&0&\U
\end{array}\right)V^*_{q-1}\cdots V^*_r.
$$

By Theorem \ref{t3.3}, $\Psi$ is determined uniquely by $\Phi$ modulo constant
unitary factors. Hence, it is sufficient to show that the indices 
$k_0,\cdots,k_{r-1}$ are uniquely determined by $\Phi$.

It follows easily from (\ref{4.6}) that 
$$
k_0=\min\left\{\kappa:~
\dim\{f\in H^2(\C^n):~\|H_{z^\kappa\Phi}f\|_2=t_0\|f\|_2\}=0\right\}.
$$
Let now $d$ be the number of indices among $k_0,\cdots,k_{r-1}$ that are to equal to
$k_0$. It follows easily from (\ref{4.6}) that
$$
d=\dim\{f\in H^2(\C^n):~\|H_{z^{k_0-1}\Phi}f\|_2=t_0\|f\|_2\}.
$$
Next, if $d<r$, then it follows from (\ref{4.6}) that
$$
k_d=\min\left\{\kappa:~
\dim\{f\in H^2(\C^n):~\|H_{z^\kappa\Phi}f\|_2=t_0\|f\|_2\}=d(k_0-\kappa)\right\}.
$$
Similarly, we can determine the multiplicity of the index $k_d$, then the next 
largest index, {\it etc}. $\bl$

\begin{cor}
\label{t4.10}
Let $\Phi\in L^\be(\mm_{m,n})$. Suppose that {\em$\|H_\Phi\|_{\text e}$}
is less than the largest nonzero superoptimal singular value of $\Phi$.
Then the indices of a monotone thematic factorization of $\Phi-\A\Phi$ are
uniquely determined by $\Phi$.
\end{cor}

\begin{cor}
\label{t4.11}
Let $\Phi\in(H^\be+C)(\mm_{m,n})$. 
Then the indices of a monotone thematic factorization of $\Phi-\A\Phi$ are
uniquely determined by $\Phi$.
\end{cor}

\

\

\noindent
Department of Mathematics
\newline
Kansas State University
\newline
Manhattan, Kansas 66506
\newline
USA

\end{document}